\documentclass[12pt,a4paper,twoside,final,notitlepage, leqno]{article} 
\usepackage[english]{babel}
\usepackage[T1]{fontenc}  
\usepackage{graphicx}
\usepackage{amsmath,amsthm,epsfig,amsfonts,bbm}   

\def \smb {{\scriptstyle \bullet }}
\def\br {\break}
 
\newcommand{\moneq}{\vspace*{-6pt} \begin{equation} \displaystyle } 
\newcommand{\moneqstar}{\vspace*{-6pt} \begin{equation*} \displaystyle } 
\newcommand{\monendstar}{\vspace*{-6pt} \end{equation*}   }
\newcommand{\monend}{\vspace*{-6pt} \end{equation}   }
\newcommand{\beq}   {\begin{equation}}
\newcommand{\enq}   {\end{equation}}
\newcommand{\be}    {\begin{enumerate}}
\newcommand{\ee}    {\end{enumerate}}

\newcommand{\Bb}

\def\N{{\rm I}\!{\rm N}}
\def\R{{\rm I}\! {\rm R}}

\def\abs#1{\mid \! #1 \! \mid }

\def\resume{\if@twocolumn
\section*{R\'esum\'e}
\else \small
\quotation{\bf \it R\'esum\'e \rule[1mm]{1.5mm}{0.2mm}\vspace{0pt}}
\fi}
\def\endresume{\if@twocolumn\else\endquotation\fi}
\def\abstract{\if@twocolumn
\noindent\section*{{\bf Abstract}}
\else \small
\quotation{\noindent \bf {Abstract.} \rule[1mm]{1.5mm}{0.2mm}\vspace{0pt}}
\fi}
\def\endabstract{\if@twocolumn\else\endquotation\fi}

\def\section*#1{}

\begin{document} 
\title{\bf \LARGE     ~ \vspace{1cm}  ~\\   Unconditionnally stable scheme   \\ ~   
  for Riccati equation ~\\~   }

\author { { \large  Fran\c{c}ois Dubois~$^{a}$ and   Abdelkader  Sa\"{\i}di~$^{b}$}     \\ ~\\ 
{\it  \small $^a$   Conservatoire National des Arts et M\'etiers}  \\
{\it  \small 15 rue Marat, F-78 210 Saint Cyr l'Ecole, France. }  \\  
{\it \small  $^b$  Institut de Recherche Math\'ematique Avanc\'ee} \\ 
{\it \small    Universit\'e Louis Pasteur } \\ 
{\it \small    7 Rue Ren\'e-Descartes, 67084 Strasbourg Cedex, France.   } \\  ~\\  ~\\   }  
\date{ July 2000~\protect\footnote{~Published, 
  {\it ESAIM: Proceedings} (http://www.esaim-proc.org), 
volume 8, ``Contr\^ole  des syst\`emes gouvern\'es
par des \'equations aux d\'eriv\'ees partielles'', 
Francis Conrad and  Marius Tucsnak Editors, p.~39-52, DOI: 10.1051/proc:2000003, 2000. 
Present   edition 18~January 2011.   } ~\\ }


\maketitle

\begin{abstract} 
In this contribution  we present a numerical scheme for the resolution
of matrix Riccati equation used in control problems. The scheme is unconditionnally
stable and the solution is definite positive at each time step of the
resolution. We prove the convergence in the scalar case and present several
numerical experiments for classical test cases.    
\end{abstract}

\bigskip \bigskip  \newpage \noindent {\bf \large 1) \quad  Introduction}

\smallskip \noindent  
We study the optimal control of a differential linear system
\moneq \label{1} 
{{{\rm d}y}\over{{\rm d}t}} \,\,=\,\, A \, y \,+\, B \, v\,\,, 
\monend
where the state variable $\, y(t) \,$ belongs to $\R^n$ and the control function 
$\, v(\smb)\,$  takes
its values in $\R^m$, with $n$ and $m$ beeing given integers. Matrix $A$ is composed by $n$
lines and $n$ columns and matrix $B$ contains $n$ lines and $m$ columns. Both matrices
$A$ and $B$ are independent of time. With the ordinary differential equation (\ref{1}) is
associated an initial condition
\moneq \label{2} 
y(0) \,=\, y_0  
\monend
with $\, y_0\,$  given in $\R^n$ and the solution of system  (\ref{1})(\ref{2}) 
 is parametrized by the function $v(\smb)$: 
The control problem consists of finding the minimum $\, u(\smb)\,$  of some
quadratic functional $\, J(\smb)$:
\moneq \label{3} 
J(u(\smb)) \,\, \leq \,\, J(v(\smb)) ,\quad \forall \, v(\smb) \,.
\monend
The functional $ \, J(\smb) \,$  depends on the control variable function
$ \, v(\smb) ,\,$   is defined by the horizon $T > 0$, 
the symmetric semi-definite positive $n$ by $n$ constant matrix $Q$
and the symmetric definite positive $m$ by $m$ constant matrix $R$. We set classically :
\moneq \label{4} 
J(v(\smb))\,\,=\,\, {1 \over 2} \int_{0}^{T} (Q y(t),y(t)) \, {\rm d}t 
\,+\,  {1 \over 2} \int_{0}^{T} (R v(t),v(t)) \, {\rm d}t  \,. 
\monend

\smallskip  \noindent $\bullet$ \quad 
 Problem  (\ref{1})(\ref{2})(\ref{3})(\ref{4})
 is a classical linear quadratic mathematical modelling of
dynamical systems in automatics (see {\it e.g.}  Lewis \cite{Le86}). 
When the control function $\, v(\smb)\,$   is supposed to be square integrable 
($v(\smb) \in {\rm L}^2 (]0; T[, \R^m) $) then the control
problem  (\ref{1})(\ref{2})(\ref{3})(\ref{4}) 
has a unique solution  $\, u(\smb) \in {\rm L}^2 (]0; T[,  \R^m) \, $
 (see for instance Lions \cite{Li68}). 
When there is no constraint on the control variable the minimum
$\, u(\smb)\,$  of the functional $\, J(v) \,$   is characterized by the condition:
\moneq \label{5}  
{\rm d}J(u) \, {\scriptstyle \bullet} \,  w \,\, = \,\,0 \, ,\quad \forall \,  w  \in  
L^2(]0,T[,\, \R^m)\,,
\monend
which is not obvious to compute directly.

\bigskip \noindent $\bullet$ \quad 
 When we introduce the differential equation  (\ref{1}) as a constraint between 
$ \, y(\smb)\,$  and $ \, v(\smb)$, 
the associated Lagrange multiplyer $ \, p(\smb)$ is a function of time and is classically
named the adjoint variable. Research of a minimum for  $ \, J(\smb)\,$  
(condition (\ref{5}))   can be
rewritten in the form of research of a saddle point and the evolution equation for
the adjoint variable is classical (see {\it e.g.}  Lewis \cite{Le86}): 
\moneq \label{6}  
{{{\rm d}p} \over {{\rm d}t}}\,+\,A^{\rm \displaystyle t} p\,+\,Q\, y\,\,=\,\,0 \,,
\monend
with a final condition at $t = T$,
\moneq \label{7}  
 p(T) \,\,=\,\,0
\monend
and the optimal control in terms of the adjoint state $\,  p(\smb) \,$  takes the form:
\moneq \label{8}  
R \, u(t) \,+\, B^{\rm \displaystyle t} \, p(t) \,\,=\,\, 0 \, .
\monend

\smallskip \noindent $\bullet$ \quad 
 We observe that the differential system  (\ref{1})(\ref{6}) together with the initial
 condition  (\ref{2}) and the final condition  (\ref{7}) 
is coupled through the optimality condition  (\ref{8}). In
practice, we need a linear feedback function of the state variable $ \, y(t) \,$ 
 instead of the adjoint variable $ \, p(t) $. 
Because adjoint state  $ \, p(\smb) \, $    depends linearily on state variable
  $ \, y(\smb) \, $  we can set: 
\moneqstar 
p(t) \,=\,  X(T-t) \, \smb \, y(t) \,, \qquad  0 \leq  t \leq  T \,, 
\monendstar
with a symmetric $n$ by $n$ matrix $ \, X(\smb)\,$ 
 which is positive definite. The final condition  (\ref{7})  is realized for each
value $ \, y(T)$, then we have the following condition:
\moneq \label{9}  
X(0) \,\,=\,\, 0 \, .
\monend
  We set $ \, K =  B R^{-1} \, B^{\rm \displaystyle t}$; 
we remark that matrix $K$ is symmetric positive definite,
we replace the control $ \, u(t) \,$  by its value obtained in relation  (\ref{8}) 
 and we deduce after  elementary algebra the evolution equation for the transition 
matrix $ \, X(\smb)$:
\moneq \label{10}  
{{{\rm d}X} \over {{\rm d}t}}\,-\, \bigl( \, XA\,+\,A^{\rm \displaystyle t}\, X \, \bigr)
\,\, +\, \,X\,K\, X\,-\,Q\,\,=\,\,0   \,, 
\monend
 which defines the Riccati equation associated with the control problem 
 (\ref{1})(\ref{2})(\ref{3})(\ref{4}).

\bigskip \noindent $\bullet$ \quad 
 In this paper we study the numerical approximation of differential system
 (\ref{9})(\ref{10}). 
Recall that the given matrices $Q$ and $K$ are $n \times  n$ symmetric matrices, with
$Q$ semi-definite positive and $K$ positive definite; the matrix $A$ is an $n$ by $n$ matrix
without any other condition and the unknown matrix $ \, X(t) \, $ 
is symmetric. We have  the following property  (see {\it e.g.}  Lewis \cite{Le86}).

\bigskip  \noindent {\bf Proposition 1. \quad The solution
of Riccati equation is positive definite.}

\noindent  
Let $K$, $Q$, $A$  be given $n \times n$  matrices with $K$,  $Q$  symmetric, 
$Q$ positive and $K$ definite positive. Let $ \, X(\smb) \,$ be the solution 
of the Riccati differential equation (\ref{10}) 
with initial condition  (\ref{9}). Then $ \,   X(t) \,$ is well defined for any 
$ \, t \geq 0, \,$ 
is symmetric and for each $ \, t > 0$, $ \, X(t) \,$ is definite positive and 
tends to a definite
positive matrix $ \, X_{\infty} \,$  as $t$ tends to infinity: 
$ \, X(t) \longrightarrow   X_{\infty} \,$  if $ \, t  \longrightarrow \infty .\,$ 
Matrix $ \,  X_{\infty} \,$  is the unique positive symmetric matrix 
which is solution of the so-called
algebraic Riccati equation:
\moneqstar  
-(XA\,+\,A^{\rm \displaystyle t}X)\,\, +\, \,XKX\,-\,Q\,\,=\,\,0\,.
\monendstar  
%

\smallskip \noindent $\bullet$ \quad 
 As a consequence of this proposition it is usefull to simplify the feedback command
law   (\ref{8}) by the associated limit command obtained by taking 
 $ \, t  \longrightarrow \infty ,\,$   that is:
\moneq \label{11} 
v(t) \,\,=\,\, -R^{-1} \, B^{\rm \displaystyle t} \, X_{\infty} \, y(t) \, ,
\monend
and the differential system   (\ref{1}) (\ref{11})
 is stable   (see {\it e.g.}  \cite{Le86}).  The practical computation
of matrix  $\, X_{\infty} \,$   by direct methods is not obvious and we refer 
{\it e.g.}  to Laub \cite{La79}. 
 If we wish to compute directly a numerical solution of instationnary Riccati
equation  (\ref{10}) classical methods for ordinary differential equations like 
{\it e.g.}  the forward Euler method
\moneqstar  
{{1} \over {\Delta t}}(X_{j+1}\,-\,X_j)\,+\,X_j \, K \, X_j\,\, -\,\,(A^{\displaystyle
t} \, X_j\, +\,X_j A)\,-\,Q\,\,=\,\,0 \,,
\monendstar 
or Runge Kutta method fail to maintain positivity of the iterate $ \, X_{j+1} \,$ 
 at the order\br 
 $(j + 1)$:
\moneq \label{12} 
(X_{j+1} \, x\,,\, x) \, > \, 0,\quad \forall \,  x \in \R^n,\quad x\,  \neq \,0 \, , 
\monend
if $ \, X_j \,$  is positive definite and if time step 
$ \, \Delta t > 0 \,$   is not small enough (see {\it e.g.}
Dieci and Eirola \cite{DE96}). Morever, there is to our best knowledge no simple way to
determine {\it a priori}  if time step  $ \, \Delta t > 0 \,$ 
is compatible or not with condition (\ref{12}).

\bigskip \noindent $\bullet$ \quad 
 In the following, we propose a method for numerical integration of Riccati
equation  (\ref{10}) which maintains condition  (\ref{12}) for each time step 
$ \, \Delta t > 0$. We
present in second section the simple case of scalar Riccati equation and present the
numerical scheme and its principal properties of the general case in section 3. We
describe several numerical experiments in section 4.

\bigskip \bigskip  \noindent {\bf \large 2) \quad   Scalar Riccati equation}

\smallskip \noindent $\bullet$ \quad 
When the unknown is a scalar variable, we write Riccati equation in the following
form:
\moneq \label{13} 
{{{\rm d}x} \over {{\rm d}t}} \,+\,k \, x^2\,-\,2\, a \, x\,-\,q\,\,=\,\,0 \,, 
\monend
with
\moneq \label{14} 
k>0,\quad q\, \geq \, 0\,, 
\monend
and an initial condition:
\moneq \label{15} 
x(0)\,\,=\,\,d,\quad d\, \geq \, 0   \,. 
\monend
We approach the ordinary differential equation (13) with a finite difference scheme of
the type proposed by Baraille \cite{Ba91} for hypersonic chemical kinetics and independently
with the ``family method''  proposed by Cariolle  \cite{Ca79}  and studied by Miellou
 \cite{Mi84}. We suppose that time step $ \,  \Delta t \, $ is strictly positive. 
The idea is to write the
approximation $ \, x_{j+1} \,$  at time step $ \, (j+1) \Delta t \, $ 
 as a rational fraction of   $ \, x_{j} \,$ with positive coefficients. 
We decompose first the real number $a$ into positive and negative parts :
$ \, a = a^+ - a ^- \,$; $ \, a^+ = \max (0; a) \geq  0 $, 
$ \,  a^-  = \max (0; - a) \geq  0, \,$ $ \,  a^+ \, a ^- =  0 \,$ 
 and factorize the product $ \, x^2 \,$  into the very simple form: 
\moneqstar  
\big( x^2 \big)_{j+1/2} \,=\, \, x_{j} \, \,  x_{j+1} \, . 
\monendstar
%

\bigskip  \noindent {\bf Definition 1. \quad  Numerical scheme in the scalar case.} 

\noindent 
For resolution of the scalar differential equation (13), we define our
numerical scheme by the following relation:
\moneq \label{16} 
{{x_{j+1}\,-\,x_{j}} \over {\Delta t}} \,+\, k \, x_{j} \,
x_{j+1}\,\,-\,\,2 \, a^+\,x_{j}\,\, +\, \,2 \,a^- \,x_{j+1}\,\, -\,\, q\,\,=\,\,0\,. 
\monend   
%

\smallskip \noindent $\bullet$ \quad 
 The scheme   (\ref{16}) is implicit because some linear equation has to be solved to
compute $ \, x_{j+1} \, $  when $ \, x_j \, $  is supposed to be given. 
In the case of our scheme this
equation is linear and the solution $ \, x_{j+1} \, $ is obtained from scheme   (\ref{16}) 
 by the homographic  relation:
\moneq \label{17} 
x_{j+1}\,\,=\,\,{{{\bigl( \, 1\,+\,2a^+\,\Delta t \, \bigr) \, x_j \,\,+\,\, q \,
\Delta t}} \over  {k\,\Delta t\,x_j \,+\,(1\,+\,2a^-\,\Delta t)}} \,. 
\monend
%

\bigskip  \noindent {\bf Proposition 2. \quad Algebraic properties of the
scalar homographic scheme.}
 
\noindent 
Let $ \, (x_j)_{j \in \N} \,$ be the sequence defined by initial condition : 
$ \, x_0 = x(0) = d \,$  and recurrence relation  (\ref{17}). 
Then sequence  $ \, (x_j)_{j \in \N} \,$ is globally defined and remains
positive for each time step: $ \, x_j > 0,$ $  \forall j \in \N, $ 
$ \, \forall \Delta t > 0 .$ If $\, \Delta t > 0 \,$ 
 is chosen such that:
\moneq \label{18} 
1\,+\,2|a|\Delta t \,-\, k\, q\, \Delta t^2 \,  \neq \, 0 \,, 
\monend
then  $ \, (x_j)_{j \in \N} \,$ converges towards the positive solution 
$x^*$  of the ``algebraic Riccati equation''   
\moneqstar
k \, x^2 \,-\, 2 \, a \, x \, - \, q \,\, = \,\, 0 
\monendstar
and 
\moneq \label{19} 
x^* \,\,=\,\, {{1} \over {k}} \Big(  a\,+\, \sqrt{a^2 \,+\, kq} \, \Big) \, . 
\monend
\smallskip \noindent $\bullet$ \quad 
 In the exceptional case where   $\, \Delta t > 0 \,$   is chosen such that
  (\ref{18})  is not satisfied, then the sequence  $ \, (x_j)_{j \in \N} \,$
is equal to the constant  
$ \, {{ 1 \, + \, 2 \, a^{+} \, \Delta t} \over  { k \, \Delta t}}\, $  for $\,j \geq  1\, $ 
and the scheme (\ref{17})  cannot be used for the approximation of Riccati equation
  (\ref{13}).  

\bigskip  \noindent {\bf Theorem 1. \quad  Convergence of the scalar scheme.} \quad 

\noindent 
We suppose that the data $\,k,a,q\, $ of Riccati equation satisfy  (\ref{14})  and
 (\ref{18}) and  that the datum $\,d\, $ of  condition
  (\ref{15})  is relatively closed to $\,x^*\, $, {\it  i.e.}:
\moneq \label{20} 
-{{1} \over {k\, \tau}} \,+\, \eta \, \,\,  \leq \, \, \,  d - x^* \, \, \,  \leq \,
\,\, C \,,
\monend
where $\,C\, $ is some given strictly positive constant $\,(C > 0)\, $, $\,x^*\, $ 
calculated according to relation   (\ref{19}) is the limit in time of the Riccati equation,
$\,\tau\, $ is defined from data $\,k, \, a, \, q\, $ by: 
\moneqstar  
\tau\,\,=\,\,{{1} \over {2\,\sqrt{a^2\,+\,kq}}} \,,
\monendstar
and $\,\eta\, $ is some constant chosen such that 
\moneq \label{21} 
0 \, < \, \eta \, < \, {{1} \over {k\,\tau}} \,.
\monend
\smallskip  \noindent  $\bullet \quad$
We denote by  $\,x(t;d)\, $  the solution of differential equation   (\ref{13}) 
with initial condition (\ref{15}). Let $\,(x_{j}(\Delta t\,;\,d_{\Delta}))_{(j\in
\N)}\, $ be the solution of the  numerical scheme defined at the relation   (\ref{17}) and
let $\,\, d_{\Delta} \,\,$ be the initial condition: 
\moneqstar   
x_0(\Delta t\,;\,d_{\Delta})\,\,=\,\,d_{\Delta} \, .
\monendstar
We suppose that the numerical initial condition $\,\,  d_{\Delta} > 0 \,\,  $
satisfies a  condition analogous to   (\ref{20}):
\moneqstar  
- {{1} \over {k\,\tau}}\,+\,\eta \, \,\,  \leq \,\,\, d_{\Delta} \, - \, x^* \, \, 
\, \leq \,\,\, C \,, 
\monendstar
with $\,C\, $ and $\,\eta > 0\, $ equal to the constant introduced in  (\ref{20})  
and satisfying  (\ref{21}). 
 
\smallskip \noindent  $\bullet \quad$
Then the  approximated value  $\,(x_j(\Delta t \,;\, d_{\Delta}))_{j \in \N}\, $ is
arbitrarily closed to the exact value  $\,x(j \Delta t \,;  \, d)\, $ for each $\,j\,
$  as  $\,\Delta t \longrightarrow 0 \, $ and $\,d_{\Delta} \longrightarrow d\, $. 
More precisely, if $\,a \,  \neq \, 0\, $ we have the  following estimate for the
error at time equal to $\,j \Delta t\, $:
\moneqstar  
\abs{ x(j \Delta t\,;\, d) \,-\, x_j(\Delta t \,;\, d_{\Delta})}  \leq 
A\,(\Delta t \,+\, \abs{ d \,-\, d_{\Delta}}),\,\,\, \forall \,  j \in \N \, , \,  0
\,<\,  \Delta t \,  \leq \,  B 
\monendstar
with constants $\,A > 0, B > 0\, $, depending on data $\,k, a, q, \eta\, $ but
independent on time step $\,\Delta t > 0\, $ and iteration $\,j\, $. 

\smallskip \noindent  $\bullet \quad$
If $\,a \,\,=\,\, 0\, $, the scheme is second order accurate in the following sense:
\moneqstar   
|x(j \Delta t;d) \,-\, x_{j}(\Delta t;d_{\Delta})|   \leq   A\,(\Delta t^2 \,+
\, \abs{d \,-\, d_{\Delta}}),\,\,\, \forall \, j \,  \in \N \, , \,  0 < \Delta t \, 
\leq \, B      
\monendstar
with constants $\,A\, $ et $\,B\, $ independent on time step $\,\Delta t\, $ and 
iteration $j$.  

\smallskip \noindent 
A direct application of the Lax theorem for numerical scheme associated to ordinary
differential equations is not straightforward because both Riccati equation and the
numerical scheme are nonlinear. Our proof is detailed in \cite{DS2k}.

\bigskip \bigskip  \noindent {\bf \large 3) \quad   Matrix Riccati equation}

\smallskip \noindent  
In order to define a numerical scheme to solve the Riccati differential equation (10)
with initial condition  (\ref{9}) we first introduce a strictly positive real number, which
is chosen positive in such a way that the real matrix $ \, [ \mu \,  {\rm I} 
- ( A + A^{\rm \displaystyle t} ) ] \,$  is definite  positive:
\moneq \label{22} 
 {{1} \over {2}}(\mu \, x \,,\, x) \,-\, (A \, x \,,\, x) > 0
\, , \quad \forall \,  x \,  \neq \, 0 \,. 
\monend
Then we introduce the definite positive matrix $M$ wich depends on $\, \mu\,$ 
 and matrix $A$: 
\moneqstar 
M \, = \, {{1}\over{2}} \, \mu \, {\rm I} \,-\, A \, . 
\monendstar
The numerical scheme is then defined by analogy with relation
(16). We have the following decomposition :
\moneq \label{23}
A \, = \, A^+ \,-\, A^- 
\monend
with  $ \, A^+ =  {{1}\over{2}} \, \mu \, {\rm I} ,\, $ 
$ \, A^- = M ,\,$ $ \, \mu > 0 , \, $ $M$ 
positive definite. Taking as an explicit
part the positive contribution $\, A^+\,$  of the decomposition   (\ref{23})  of matrix A and in the
implicit part the negative contribution $ \, A^- = M \,$ of the decomposition 
 (\ref{23}), we get 
\moneq \label{24} \left\{ \begin{array}{l} \displaystyle \qquad 
{{1} \over {\Delta t}} (X_{j+1} \,-\, X_j) \,+\, {{1} \over {2}} (X_j K X_{j+1}
\,+\, X_{j+1} K X_j) \,+\,\\
\displaystyle  \qquad \qquad \qquad \qquad\qquad \qquad\qquad 
\,+\, (M^{\rm \displaystyle t} X_{j+1} \,+\, X_{j+1} M) \,\,=\,\, \mu
X_j \,+\, Q  \,.
 \end{array} \right. \monend 
The numerical solution given by the scheme $\, X_{j+1} \,$  
at time step $ \, j +1 \,$ is then defined
as a solution of Lyapunov matrix equation with matrix $X$ as unknown:
\moneqstar
S^{\rm \displaystyle t}_j \, X \,\,+\, \, X \, S_j \,\,=\,\, Y_j \,  
\monendstar
with 
\moneq \label{25}
S_j \, \,\,=\,\, \, {{1} \over {2}}  I \,+\, {{\Delta t} \over {2}} K X_j \,+\,
\Delta t \, M \,
\monend 
and
\moneq \label{26}
Y_j \, \,\,=\,\, \, X_j \,+\, \mu \, \Delta t \, X_j \,+\, \Delta t \, Q  \,.
\monend 
We notice that $ \, S_j \,$  
 is a (non necessarily symmetric) positive matrix and that $ \, Y_j \,$  is
a symmetric definite positive matrix if it is the case for $ \, X_j$.

\bigskip  \noindent {\bf Definition 2. \quad Symmetric matrices.} 

\noindent   
Let $n$ be an integer greater or equal to $1$. We define by $\,{\cal{S}}_n(\R)$, 
(respectively $\,{\cal{S}}^+_n(\R)$, $\,{\cal{S}}^{+*}_n(\R)$) the linear
space (respectively the closed cone, the open cone) of symmetric-matrices (respectively
symmetric positive and symmetric definite positive matrices). 
The following inclusions  $\,\, {\cal{S}}^{+*}_n(\R) \, \subset \, {\cal{S}}^+_n(\R) 
\, \subset \, {\cal{S}}_n(\R) \,\,\, $ are natural.  

\bigskip  \noindent {\bf Proposition 3.  \quad Property of the Lyapunov equation.} 

\noindent   
Let $\,S$ be a matrix which is not necessary symmetric, such that the associated
quadratic form:  $\,\R^n \ni x \longmapsto (x,Sx) \in \R $, is strictly positive 
{\it i.e.} 
\moneqstar 
S \,+\, S^{\rm \displaystyle t} \,\, \in {\cal{S}}^{+*}_n(\R) \,. 
\monendstar
Then the application $\,\, \varphi  \,\, $ defined by :
\moneq  \label{27} 
{\cal{S}}_n(\R) \, \ni \, X \, \longmapsto \, \varphi (X) \, \,\,=\,\, \,
S^{\rm \displaystyle t}\, X \,\,+\, \, X \, S \, \,\, \in \, {\cal{S}}_n(\R) \,,
\monend 
is a one to one bijective application on the space $\,{\cal{S}}_n(\R)$ of real
 symmetric matrices of order $\,n$.
Morever, if matrix $\,\varphi (X)\, $   is positive (respectively definite 
positive)  then the matrix $\,X\, $ is also  positive (respectively definite
positive): 
\moneqstar 
{\rm if}  \,\,\, \varphi (X) \in {\cal{S}}^{+}_n(\R)  \,, \quad 
{\rm then}   \,\,\,   X \in {\cal{S}}^{+}_n(\R) \, .   
\monendstar  
%

\smallskip \noindent $\bullet$ \quad  
The numerical scheme has been written as an equation with unknown 
$ \, X = X_{j+1} \, $  which takes the form: $ \, \varphi_j(X) = Y_j \,$ 
 with   $ \, \varphi_j \, $  given by a relation of the type (\ref{27}) 
 with the help of matrix $ \, S_j \, $  defined in  (\ref{25})  and a datum matrix 
$ \, Y_j \,$   defined by relation (\ref{26}). Then we have the following propositions. 

\bigskip  \noindent {\bf Proposition 4.  \quad Homographic scheme computes  
a definite  positive matrix.} 

\noindent  
The matrix $\,X_{j}\,$ defined by  numerical scheme  (\ref{24})   with the initial
condition $ \,  X_0  =  0 \,$ 
is  positive for each time step $\,\Delta t > 0  \,$ : 
\moneqstar
 X_j \in {\cal{S}}^+_n(\R), \,\,\qquad \quad  \forall \,  j \ge 1 \,  . \, 
\monendstar
If there exists some integer $ \,  m \, $ such that $\,  X_m \,$ belongs to the open
cone  $\,  {\cal{S}}^{+*}_n(\R)\,$, then  matrix $\,  X_{m+j} \,$ belongs to the open
cone  $\,  {\cal{S}}^{+*}_n(\R)\,$ for each   $\,j $.   

\bigskip  \noindent {\bf Proposition 5.  \quad Monotonicity.} 

\noindent   
Under the   condition  
\moneqstar 
{{1} \over {2} }\, \bigl( K X_{\infty} + X_{\infty} K \bigr)  \,
< \, \bigl( \mu \,+\, {{1} \over {\Delta t}}  \bigr) \, I \,, \,  
\monendstar

\smallskip \noindent  
the scheme  (\ref{24})  is monotone and we have  more precisely : 
\moneq  \label{28} 
\Bigl( \, 0 \, \,  \leq \, \, X_j \, \,  \leq \, \, X_{\infty} \,\Bigr) \,
\Longrightarrow  \, \Bigl( \, 0 \, \,   \leq \, \, X_j \, \,  \leq \, \, X_{j+1} \, \, 
\leq \, \, X_{\infty} \Bigr) \, . 
\monend    

\bigskip \bigskip  \newpage  \noindent {\bf \large 4) \quad  First numerical experiments}

\bigskip     \noindent {\bf \large 4-1   \quad  Square root function}  

\smallskip \noindent $\bullet$ \quad   
The first example studied is the resolution of the equation :
\moneq  \label{29} 
{{{\rm d}X} \over {{\rm d}t}} \,+\, X^2 \,-\, Q \,\,=\,\, 0 , \quad X(0)
\,\,=\,\, 0 \, 
\monend 
with $ \, n = 2 ,\,\, A = 0 ,\,\,  K = I \,\,$ and  matrix $\, Q \,$ equal to 
\moneq  \label{30} 
Q \,\,=\,\, {1\over2} \,\, 
 \left(\begin{array}{cc}  1 & -1 \cr   1 & 1 \end{array}\right)   \,\, 
 \left(\begin{array}{cc}  1 &  0 \cr 0 & 100  \end{array}\right)  \,\,  
 \left(\begin{array}{cc}  1 &  1 \cr   -1 & 1 \end{array}\right) \, .  
\monend 

\smallskip \noindent $\bullet$ \quad   
We have tested our numerical scheme for fixed value  $\,
\Delta t = 1/100 \,$ and   different values of parameter $\, \mu \, : \,$
$\,  \mu = 0.1 , \, 10^{-6} ,\,  10^{+6} . \,$ For small values of parameter $ \, \mu
,\,$ the behaviour of the scheme does not change between $\,   \mu = 0.1 \,$ and
$\,   \mu = 10^{-6} . \,$ Figures 1 to 4 show the evolution with time of the
eigenvalues of matrix $\, X_{j} \,$ and the convergence is achieved to the square root
of matrix $\, Q .\,$ For large value of parameter $\, \mu  \, \, ( \mu =  10^{+6} ),
\,$ we loose completely consistency of the scheme (see figures 5 and 6).

\centerline {
{\includegraphics[width=.49 \textwidth]{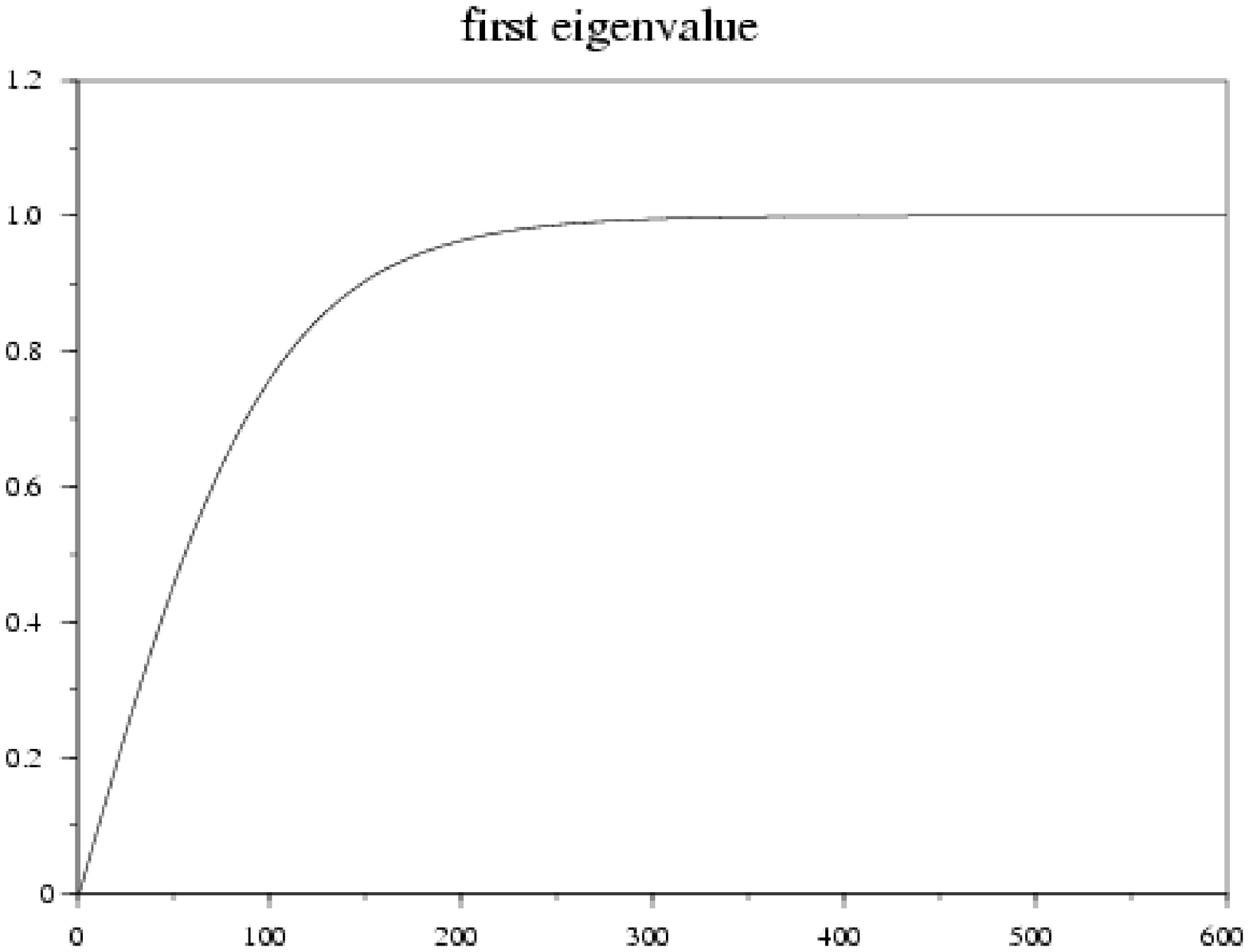}}
{\includegraphics[width=.49 \textwidth]{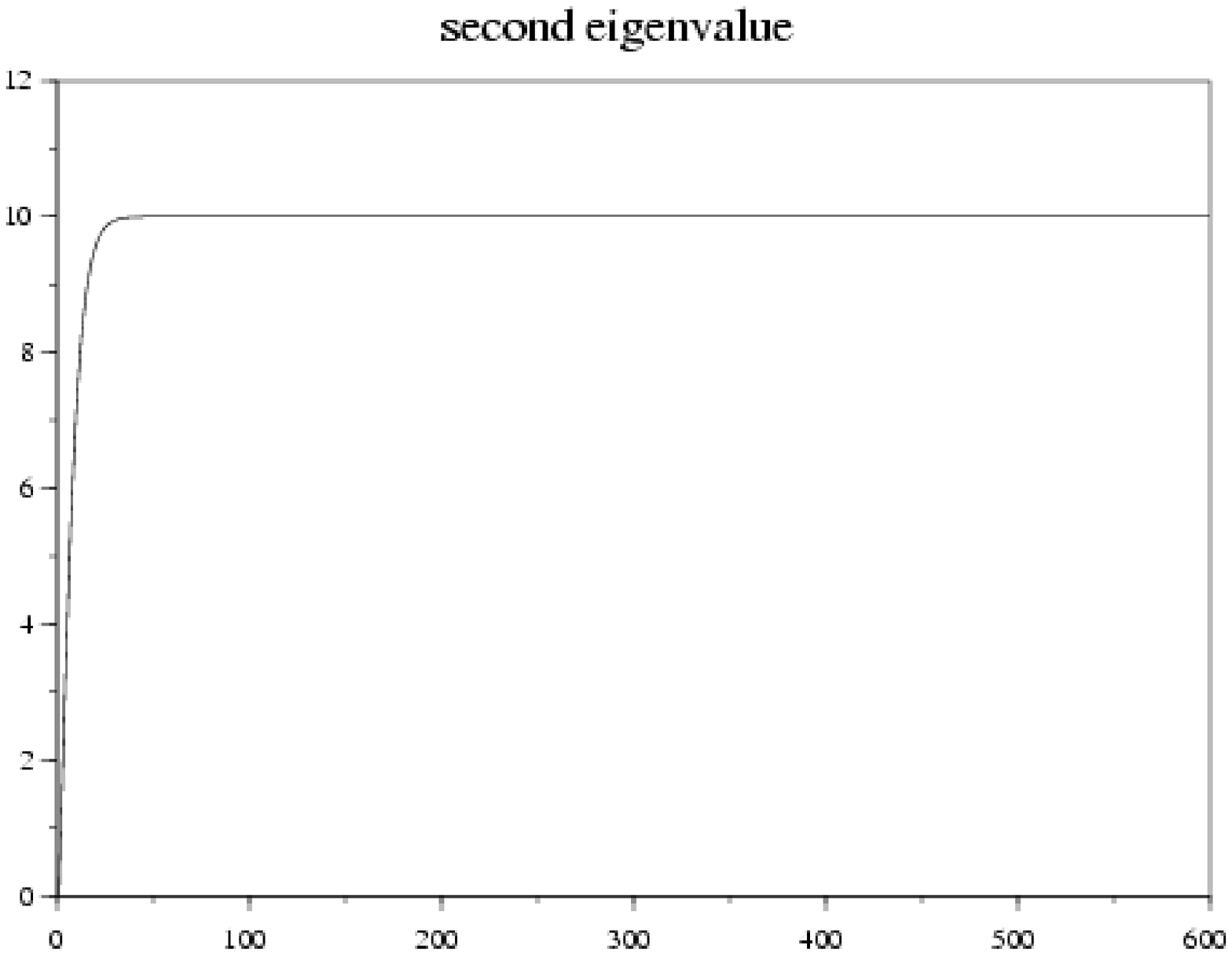}}} 
\centerline {\bf Figures 1 and 2. \quad Square root function test.}  
\centerline {Two first eigenvalues of numerical solution  ($\mu \,=\,0.1$).}     
\bigskip  

\bigskip   \centerline { 
{\includegraphics[width=.49 \textwidth]{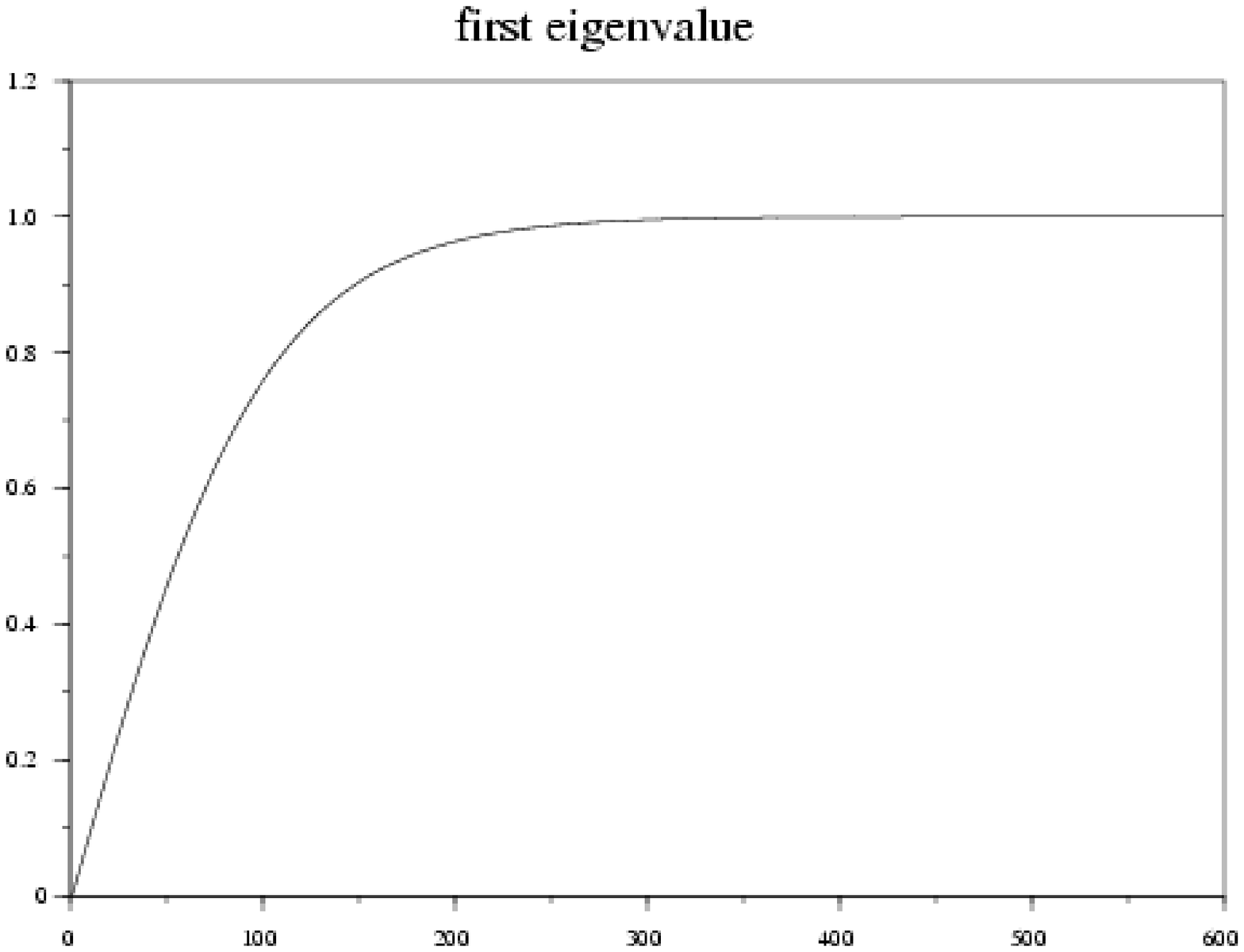}}
{\includegraphics[width=.49 \textwidth]{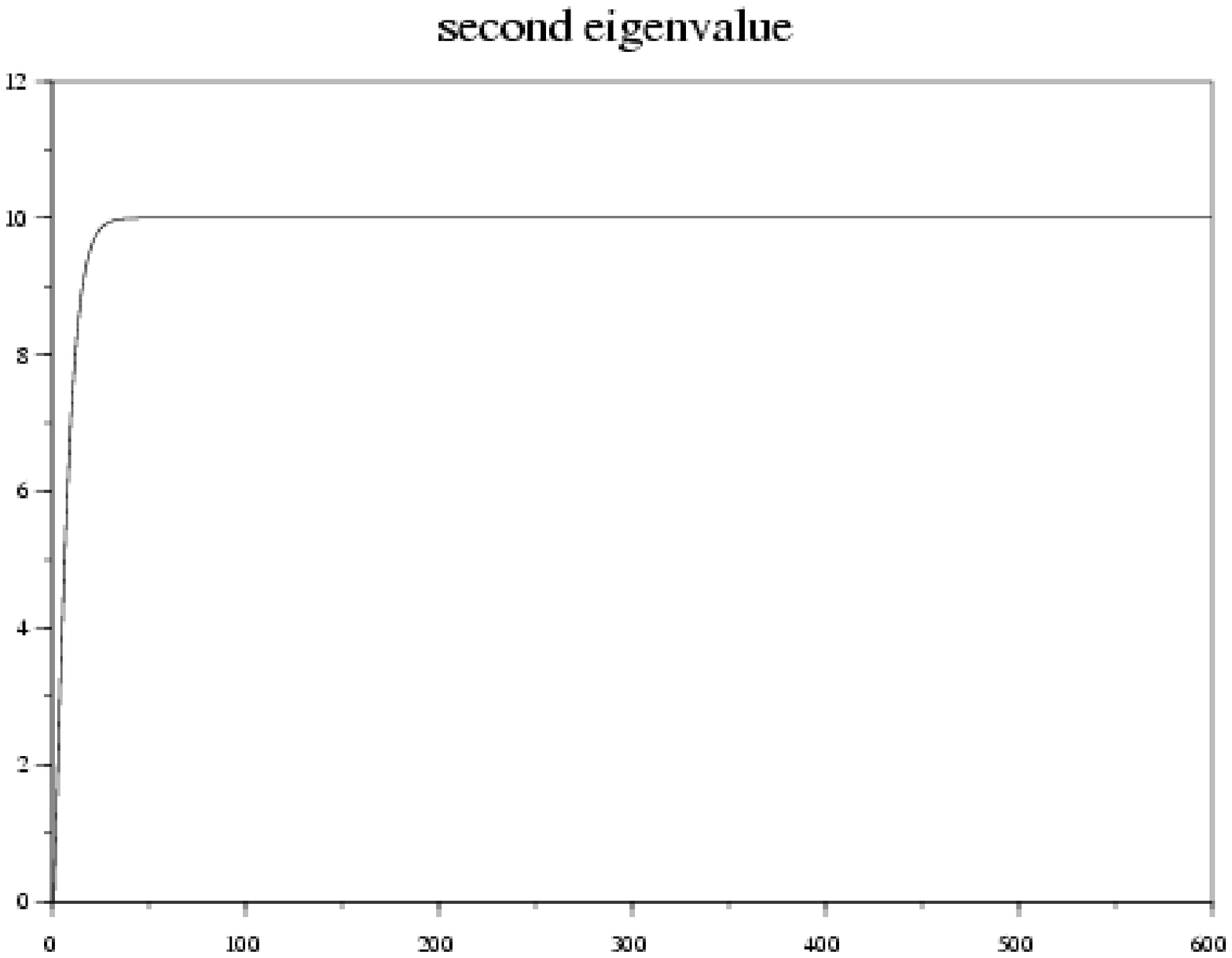}}} 
\centerline {\bf Figures 3 and 4. \quad Square root function test.}  
\centerline {Two first eigenvalues of numerical solution  ($\mu \,=\,10^{-6}$).}  
\bigskip  

\bigskip \bigskip  \newpage \centerline { 
{\includegraphics[width=.49 \textwidth]{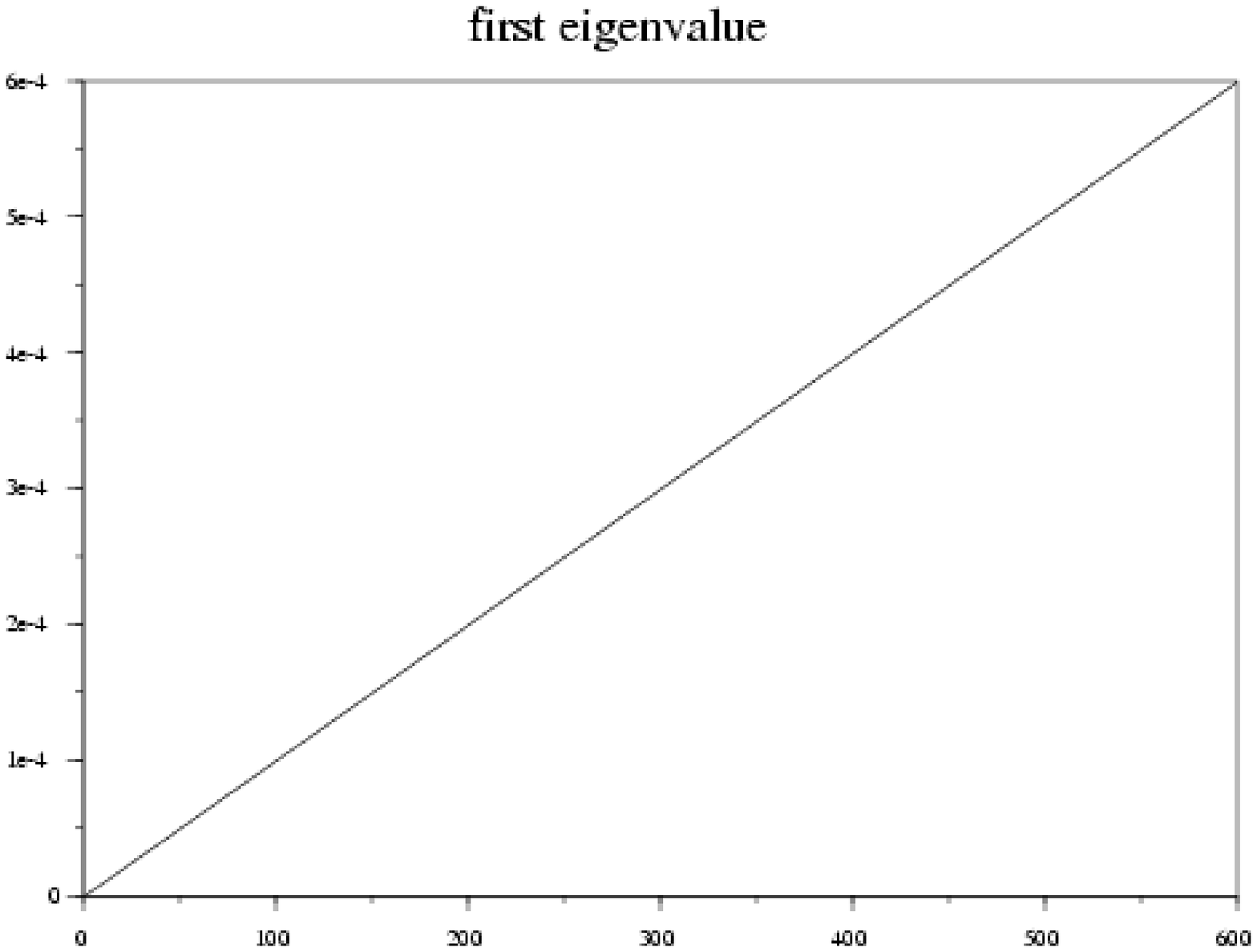}}
{\includegraphics[width=.49 \textwidth]{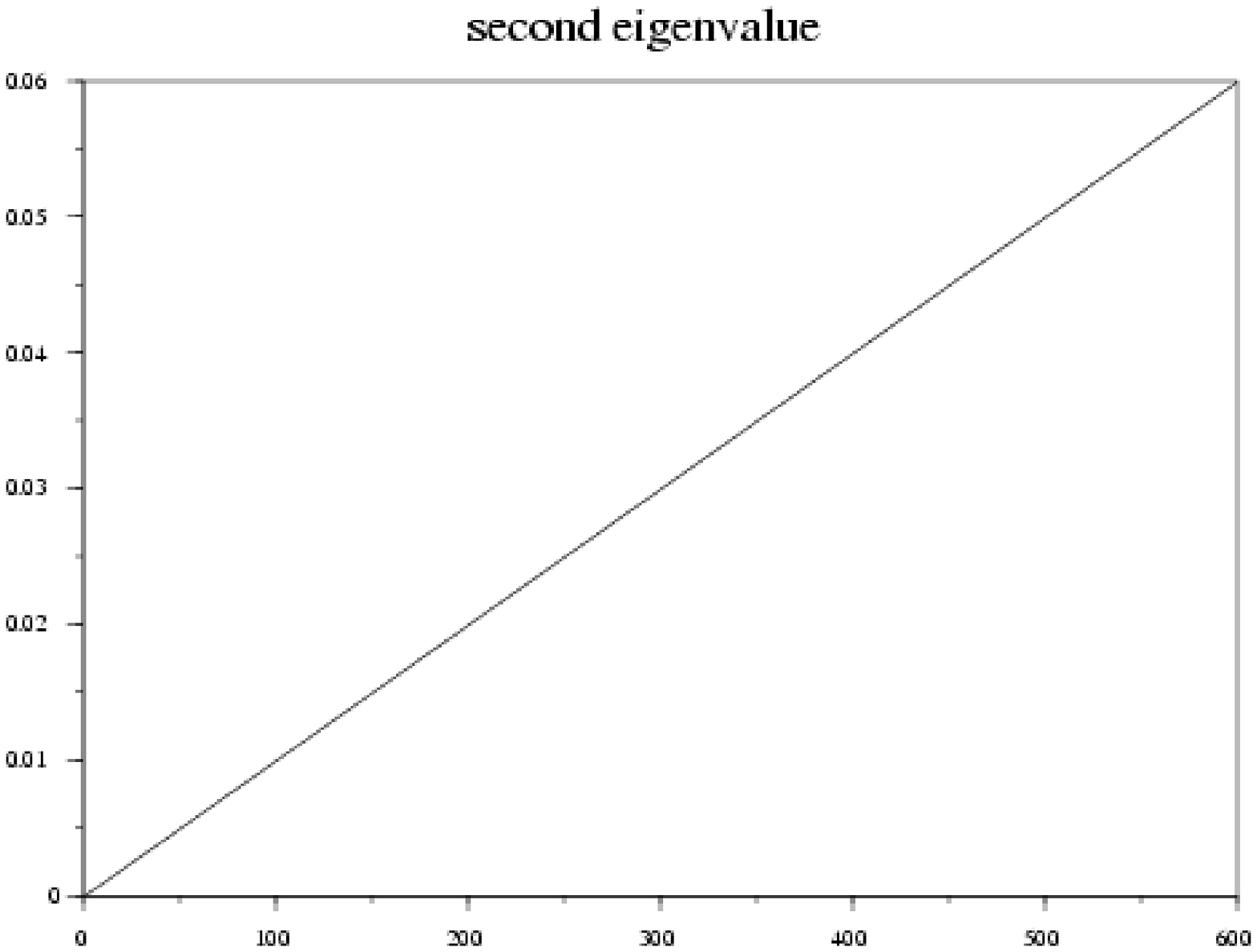}}} 
\centerline {\bf Figures 5 and 6. \quad Square root function test.}  
\centerline {Two first eigenvalues of numerical solution ($\mu \,=\,10^{+6}$).} 
\bigskip  

\bigskip     \noindent {\bf \large 4-2   \quad  Harmonic oscillator}  

\smallskip \noindent $\bullet$ \quad   
The second exemple is the classical harmonic oscillator. Dynamical system
$\, y(t) \,$ is governed by the second order differential equation with  command $\,
v(t) \,$ :
\moneq  \label{31} 
{{{\rm d}^2 y(t)} \over {{\rm d}t^2}} \,+ \, 2 \, \delta \,
{{{\rm d} y(t)} \over {{\rm d}t}} \, + \, \omega^2 \, y(t) \,\, = \,\, b \,
v(t) \,.    
\monend 
This equation is written as a first order system of differential equations :
\moneq  \label{32} 
Y \,= \,  \left(\begin{array}{c}  y(t) \cr  {\displaystyle{{\rm d} y(t)} \over
{\displaystyle {\rm d}t}}  \end{array}\right)   \,, \quad 
 { {{ \rm d} Y} \over { {\rm d}t}} \,\,= \,\,  
 \left(\begin{array}{cc}  0 & 1 \cr   - \omega^2 & -2 \, \delta  \end{array}\right) 
 \,Y(t) \,\, + \,\,   \left(\begin{array}{c}  0 \cr b \, v(t)  \end{array}\right) \,.   
\monend 
In this case, we have tested the stability of the scheme for fixed value of
parameter $\, \mu \,( \mu = 0.1) \,$ and  different values of time step 
$\, \Delta t \,$   and coefficients of matrix $R$ inside
the cost function of relation  (\ref{4}): 
\moneqstar
R \,=\,  \left(\begin{array}{cc} \alpha & 0 \cr 0 & \alpha  \end{array}\right) \, .
\monendstar

\smallskip \noindent $\bullet \quad$
We have chosen three sets of parameters : $ \, \alpha = \Delta t = 1/100 $
(reference experiment, figures 7 and 8),  $ \, \alpha = 10^{-6} \,, \,  \Delta t =
1/100 $ (very small value for $\alpha$, figures 9 and 10)  and  $ \, \alpha = 1/100
\,, \,  \Delta t = 100 $ (too large value for time step,  figures 11 and 12). Note
that for the last set of parameters, classical explicit schemes fail to give any
answer.  As in previous test case, we have represented the two eigenvalues of discrete
matrix solution $\, X_j \,$  as time is increasing. On reference experiment (figures~7
and 8), we  have convergence of the solution to the solution of algebraic Riccati
equation. If control parameter $\, \alpha \,$ is chosen too small, the
first eigenvalue of Riccati matrix oscillates during the first time steps but reach
finally the correct values of limit matrix, the solution of  algebraic Riccati
equation. If time step is too large, we still have stability but we loose also
monotonicity. Nevertheless, convergence is achieved as in previous case.

\bigskip  \newpage  \centerline { 
{\includegraphics[width=.49 \textwidth]{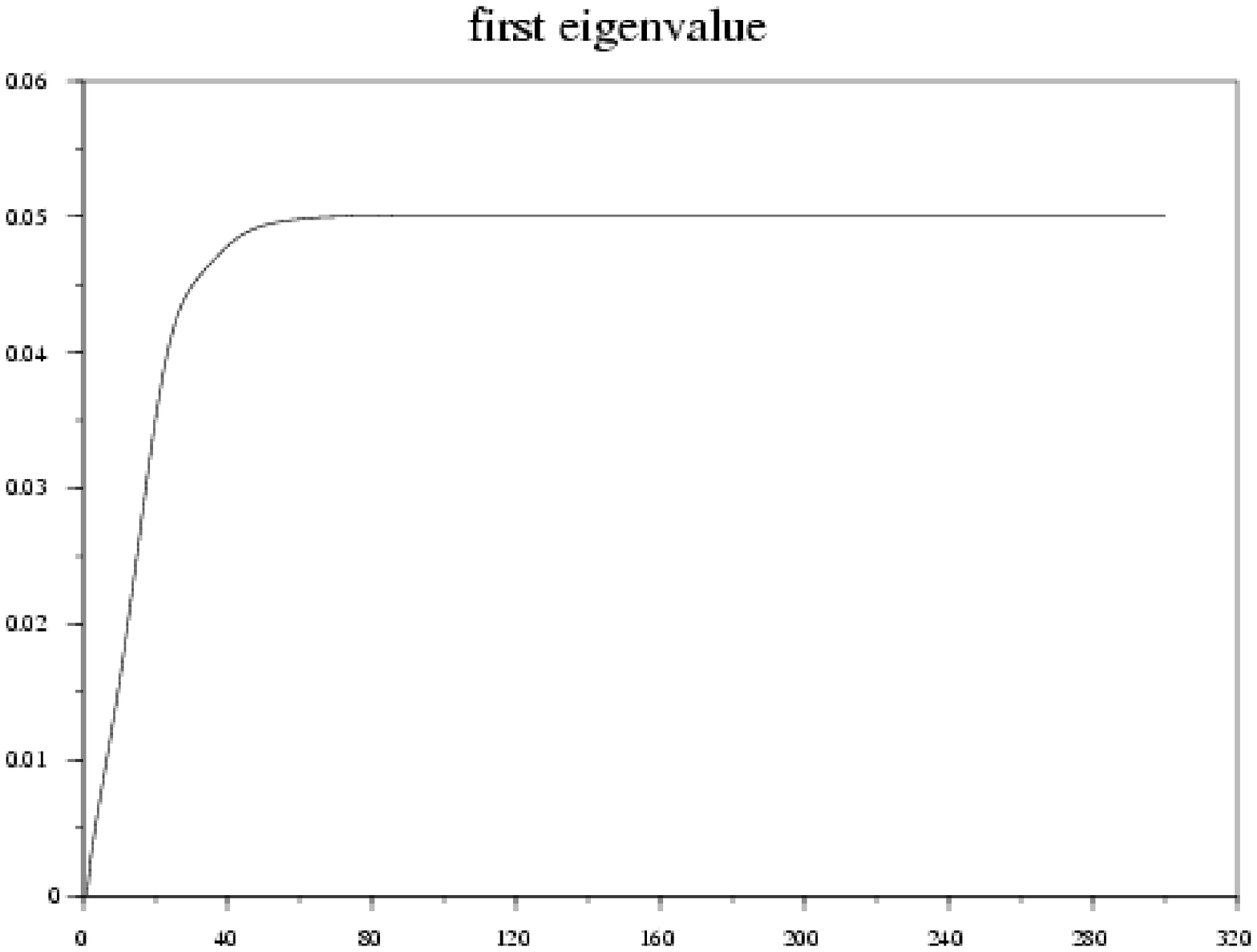}}
{\includegraphics[width=.49 \textwidth]{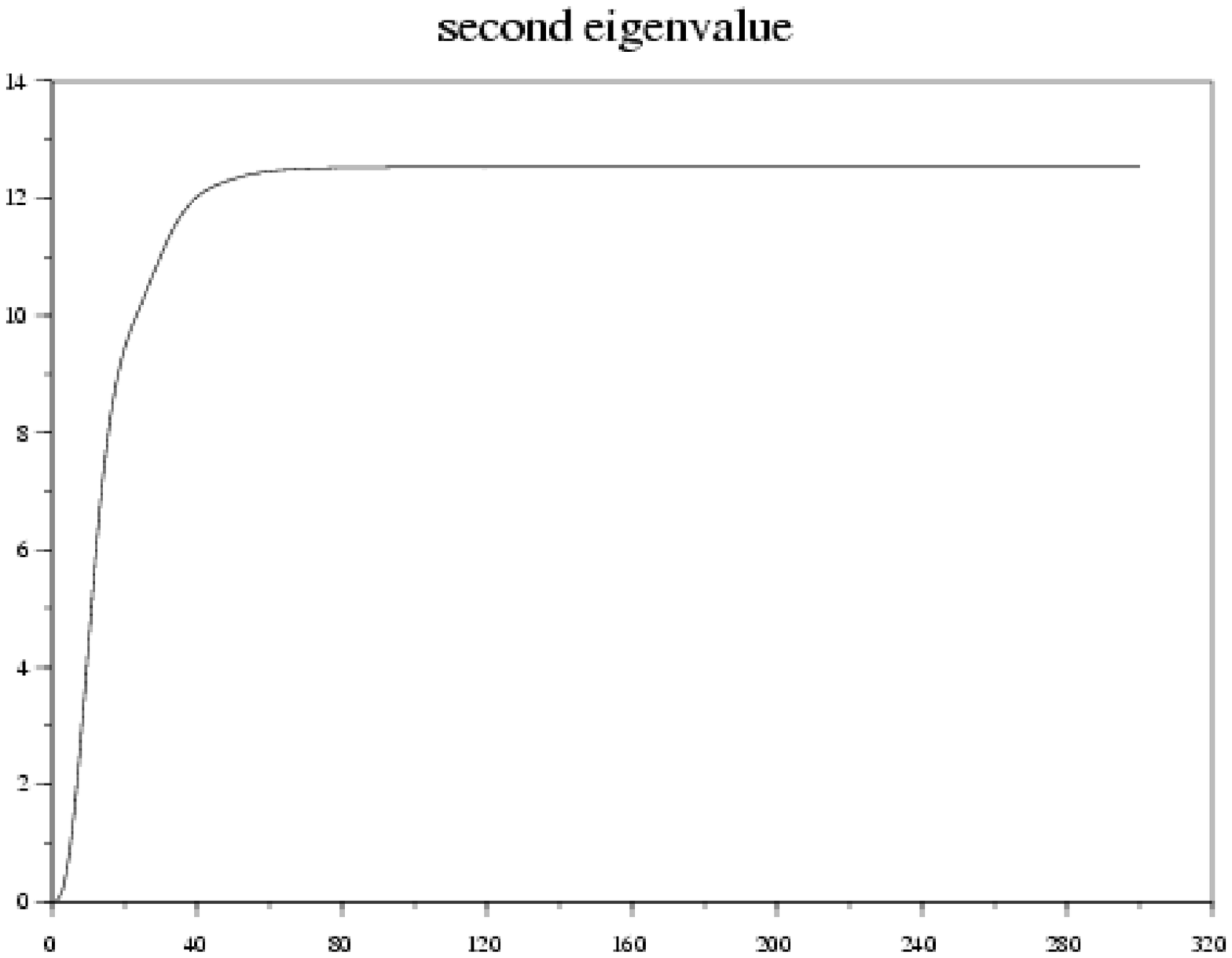}}} 
\centerline {\bf Figures 7 and 8. \quad Harmonic oscillator.} 
\centerline {Two first eigenvalues of numerical solution 
($\mu \,=\, 0.1 , \,\, \alpha \,=\,  0.01 , \,\,  \Delta t \,=\,  0.01 $).}  
\bigskip  

\bigskip  \bigskip  \centerline { 
{\includegraphics[width=.49 \textwidth]{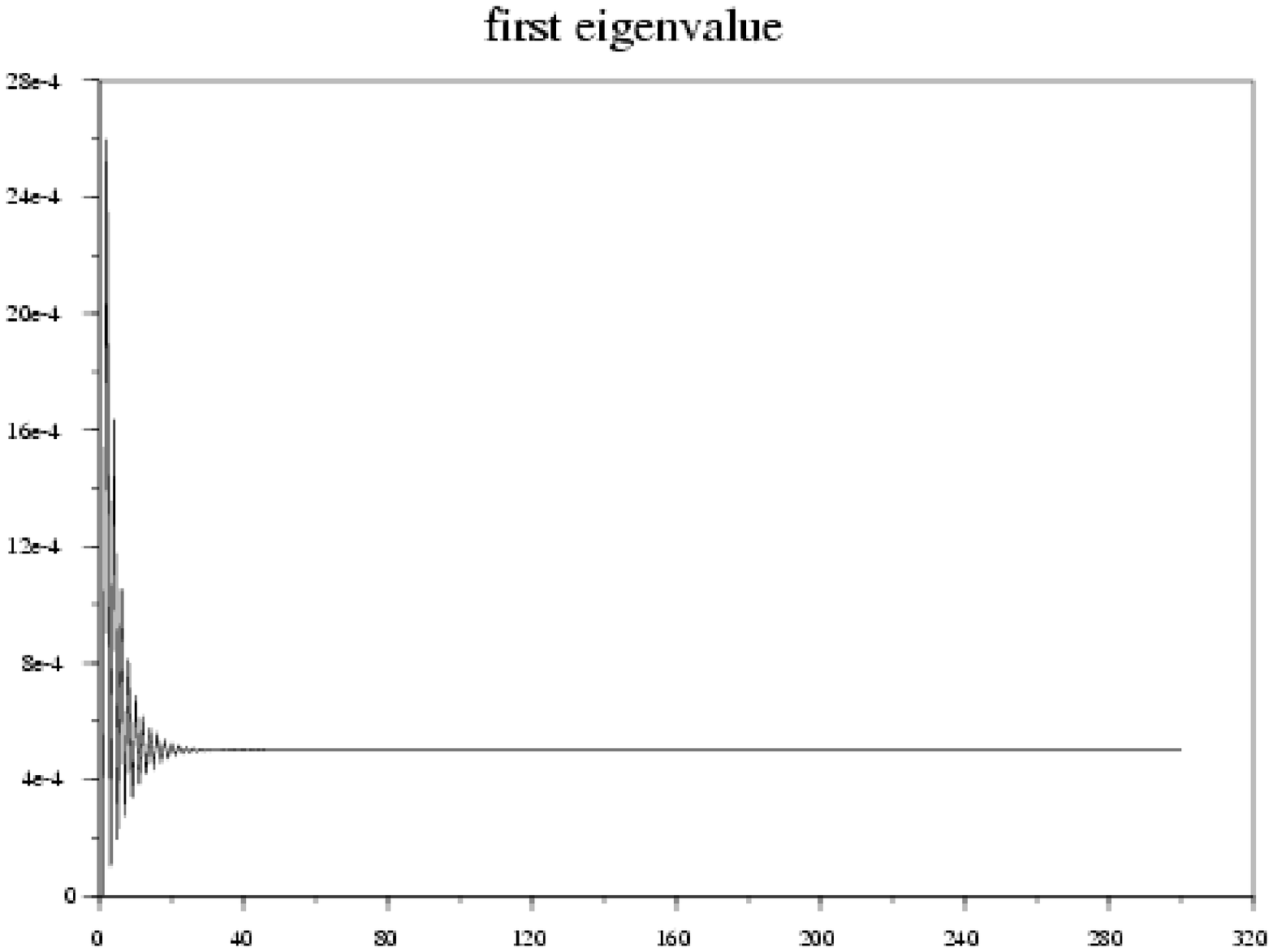}}
{\includegraphics[width=.49 \textwidth]{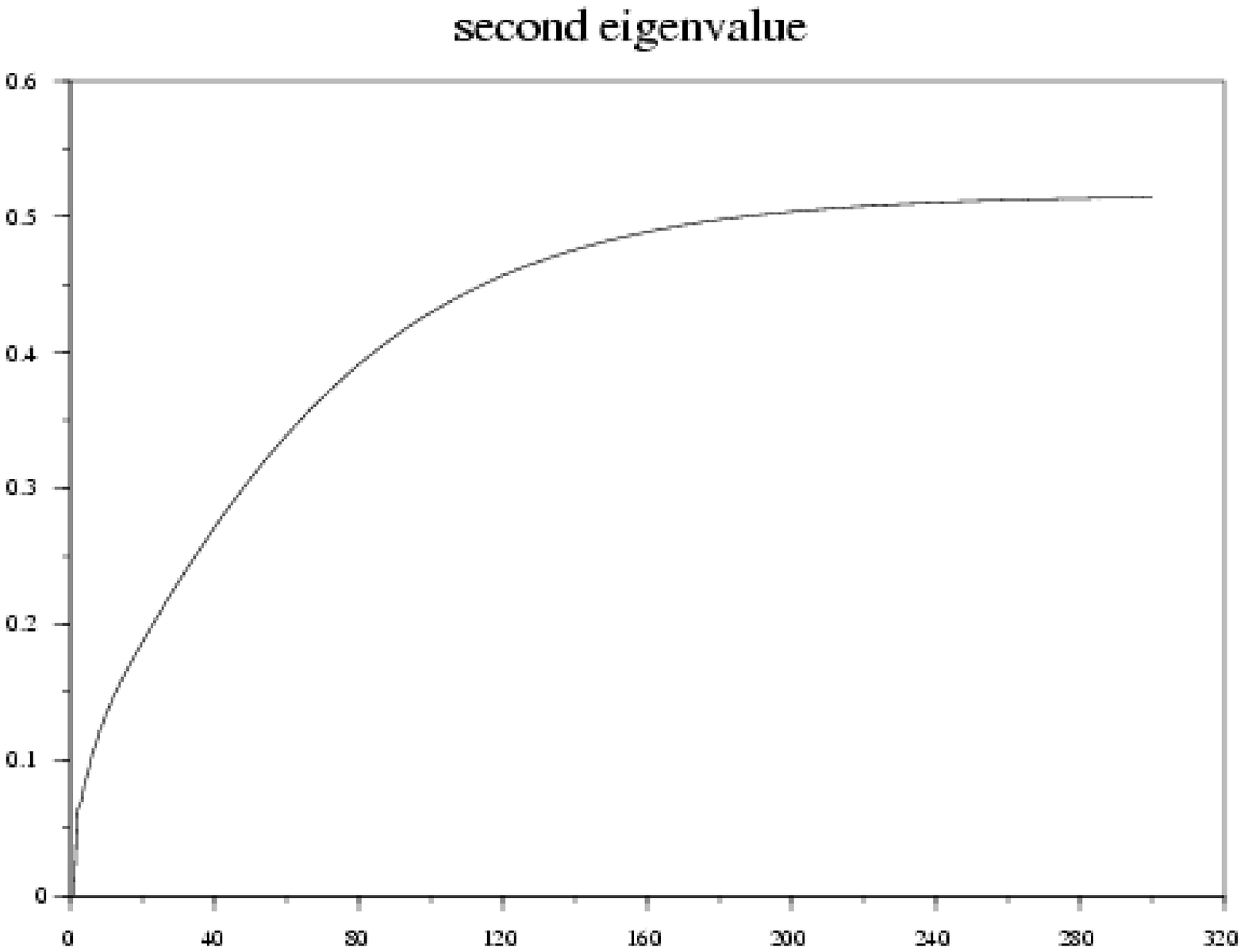}}} 
\centerline {\bf Figures 9 and 10. \quad Harmonic oscillator.} 
\centerline {Two first eigenvalues of numerical solution 
($\mu \,=\, 0.1 , \,\, \alpha \,=\, 10^{-6} , \,\,  \Delta t \,=\,  0.01 $). }
\bigskip  

\bigskip   \bigskip  \centerline { 
{\includegraphics[width=.49 \textwidth]{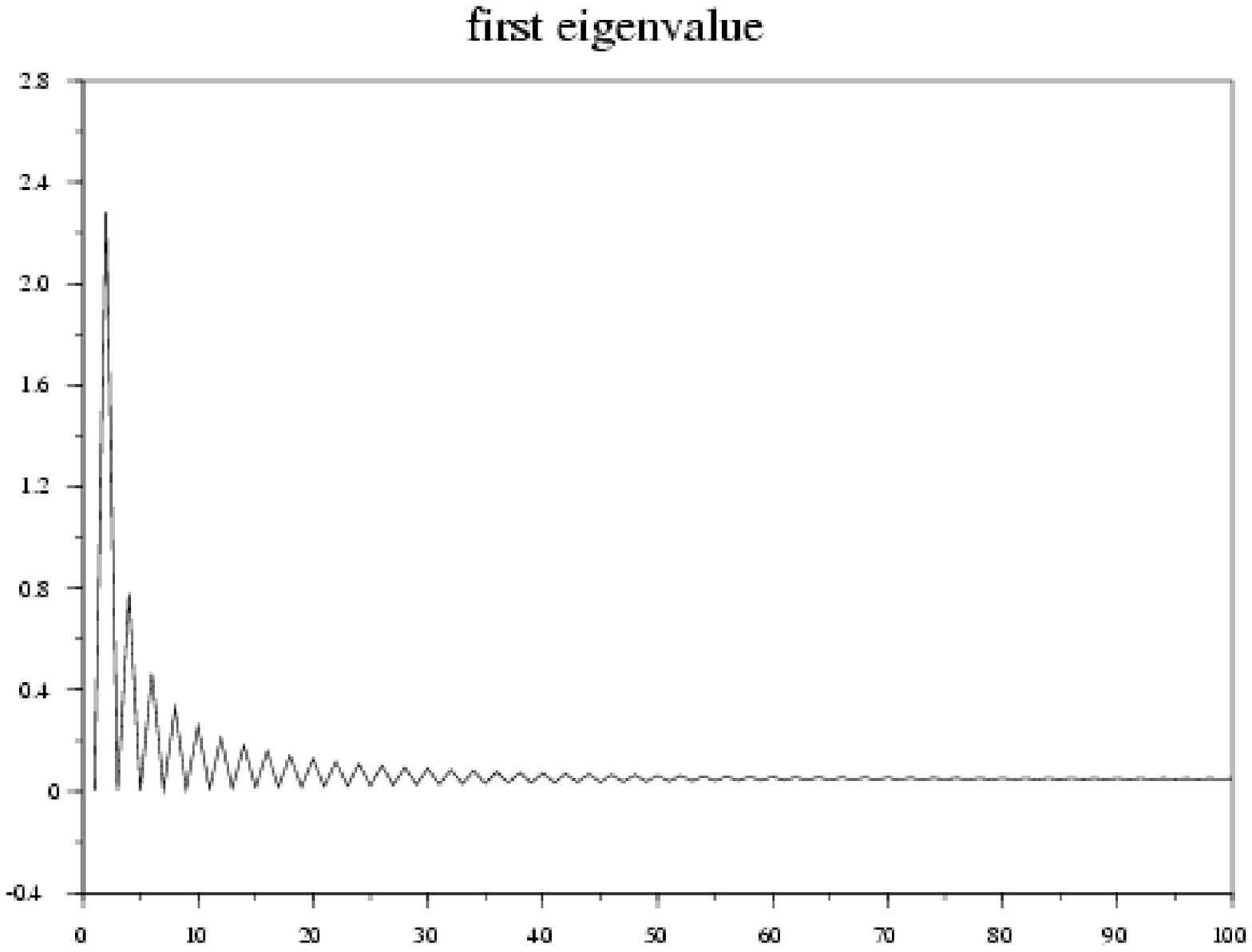}}
{\includegraphics[width=.49 \textwidth]{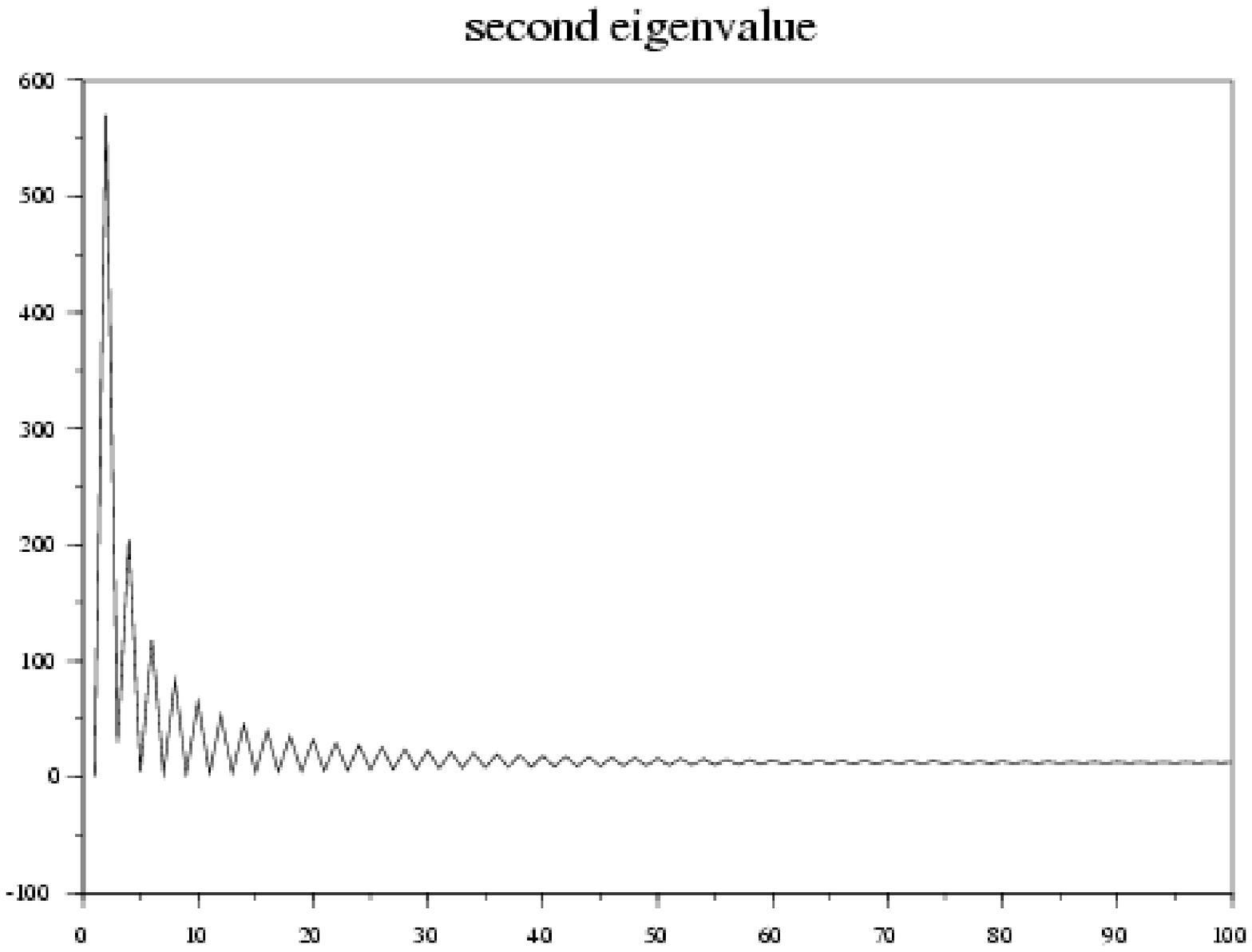}}} 
\centerline {\bf Figures 11 and 12. \quad Harmonic oscillator.} 
\centerline {Two first eigenvalues of numerical solution 
($\mu \,=\, 0.1 , \,\, \alpha \,=\,  0.01 , \,\,  \Delta t \,=\,  100 $). } 
\bigskip  

\bigskip \bigskip  \newpage \noindent {\bf \large 5) \quad   Conclusion}
  
\smallskip \noindent  
We have proposed a numerical scheme for the resolution of the matrix Riccati equation.
The scheme is implicit, unconditionnaly stable, needs to use only one scalar
parameter and to solve a linear system of equations for each time step. This scheme
is convergent in the scalar case and has good monotonicity properties in the matrix
case. Our first numerical experiments show stability and robustness when various
parameters have large variations. Situations where classical explicit schemes fail to
give a solution compatible with the property that solution of Riccati equation is a
definite positive matrix have been computed. We expect to prove convergence in
the matrix case and we will present in \cite{DS2k} experiments on realistic test models
such as a string of vehicles and the discretized wave equation.

\bigskip \bigskip  \noindent {\bf \large   Acknowledgments}

\smallskip \noindent 
The authors thank Marius Tucsnak for helpfull comments on the first draft of this
paper.


\bigskip \bigskip  

\end{document}